\documentclass[12pt]{article}
\usepackage[cp1251]{inputenc}
\usepackage{epic}
\usepackage{eepic}
\usepackage{graphics}
\usepackage{amsfonts}
  \usepackage{amsthm}
 \usepackage{amsmath}
 \usepackage{amscd}
 \usepackage{amssymb}
 \usepackage{latexsym}
\usepackage{graphicx}

 \numberwithin{equation}{section}
 \newtheorem{thm}{Theorem}
 
 \newtheorem{prop}{Proposition}
 \theoremstyle{definition}
 
\newcommand{\ds}{\displaystyle}

\begin{document}

\title{Non-Integrability of Geodesic dynamics of Chazy-Curzon space-time}

\author{Georgi Georgiev }
\date{}

\maketitle

\begin{abstract}
We study  the   integrability of the geodesic equations of the
Chazy- Curzon space-time. It was established that for the
equilibrium point $p_{\rho}=p_z=z=0$ and,  $\rho_0 \in (1,\, 2)$,
there are only periodic solutions,  the Hamiltonian system,
describing geodesic motion of Chazy-Curzon space-time has no
additional analytic first integral. Our approach is based on the
following: if the system has a family of periodic solutions around
an equilibrium and if the period function is infinitely branched
then the system has no additional analytical first integral.
\end{abstract}

\bigskip
\section{Introduction}
Let $H$ be a smooth real-valued function in $2n$ real variables
$(p,\, q),\, p,\,q\, \in {\bf R^n}$. Assume that $dH(a)=0$, where,
$a$ is an equilibrium point for the Hamiltonian system  $X_H$ (with
$n$ degrees of freedom), given by

\begin{equation}
\frac{dq}{dt}=\frac{\partial H}{\partial p},\,\,
\frac{dp}{dt}=-\frac{\partial H}{\partial q}.\nonumber
\end{equation}
We  write the Hamiltonian systems in the form

\begin{equation}
\dot{x}=X_{H}(x),\, x\in {\bf R^{2n}},\nonumber
\end{equation}
where $X_{H}$ is the hamiltonian vector field. The system is called
Liouville - Arnold  Integrable near $a$ if there exist $n$ functions
in involution $f_1=H,\, f_2,\, \dots f_n$, defined around $a$, which
are functionally independent. The Poisson bracket of $f$ and $g$ are
 \begin{equation}
\{f,\,g\}=X_{f} (g) =\sum\frac{\partial f}{\partial p}\frac{\partial
g}{\partial q}-\frac{\partial f}{\partial q}\frac{\partial
g}{\partial p}=-\{g,\, f\}.\nonumber
\end{equation}
We say that the functions $f$ and $g$ are in involution if the
Poisson bracket is equal to 0. This means that $df_1,\, df_2,\,
\dots df_n$ are linearly independent around the equilibrium $a$ and
$f_j=const$ for all $j$ define a smooth submanifold, this level
manifold is invariant under $X_{f_i}$. We have $X_{f_i}f_{j}=0$ and
$[X_{f_j},\,X_{f_k}]=X_{\{f_{j},f_{k}\}}=0$ - the vector fields
commute. The compact and connected component of $M_c:=\{
f_{j}=c_j,\, j=1,\, \dots , n\}$ is diffeomorphic to a torus.

Given a Riemannian manifold $M$, a geodesic may be defined as curve
that results from the application on the principle of least action
(we follow\cite{Dub}). A differential equation describing their
shape may be derived using variational principles, by finding
extremum of the energy of a curve. Given a smooth curve
$$\gamma \,  :\, I\rightarrow M $$
this maps an interval $I$ in real line to the $M$, one writes the
energy
$$E(\gamma )=\frac{1}{2}\int_I g(\dot{\gamma (t)}, \dot{\gamma (t)})dt,$$
where $\dot{\gamma (t)}$ is the tangent vector to the curve $\gamma$
at the point $t\in I$. Here $g(., .)$ is the metric tensor on the
manifold $M$.

In terms of local coordinates on $M$, the geodesic equation is
$$\ds{\frac{d^2 x^a}{dt^2}}+\Gamma_{bc}^a\ds{\frac{dx^b}{dt}.\frac{dx^c}{dt}}=0,$$
where the $x^a(t)$ are coordinates of $\gamma (t)$, $\Gamma_{bc}^a$
are the Christoffel symbols and indexes up and below means summation
convention.

Geodesics can be understood to be the Hamiltonian flows of a special
Hamiltonian vector field defined on the cotangent space of the
manifold. The Hamiltonian is constructed from the metric on the
manifold, and is thus a quadratic form consisting entirely of the
kinetic term. The geodesic equations are second order differential
equations, they can be expressed as a first order equations by
additional independent variables. Note that $x_a$ is the coordinates
on the neighborhood $U$ induces local trivialization of
$$T^*M|_U  \simeq U\times \bf{R^n}$$
by the map which  sends a point $\eta \in T^*_x M|_U$ of the form
$\eta=p_adx^a$ to the point $(x,\, p_x )\in U\times \bf{R^n}$. The
Hamiltonian is
$$H(x,\, p_x)=\ds{\frac{1}{2}}g^{ab}(x)p_a p_b.$$
Here $g^{ab}(x)$ is the inverse of the metric tensor
$g^{ab}(x)g_{bc}(x)=\delta^a_c$. The behavior of the metric tensor
under transformation of  coordinates implies that $H$ is invariant
under change of variables. The geodesic equations may be written as
$$\dot{x^a}=\ds{\frac{\partial H}{\partial p_a}}=g^{ab}(x)p_b,\, \dot{p_a}=-\ds{\frac{\partial H}{\partial x^a}}=-\ds{\frac{1}{2}\frac{\partial g^{ab}(x)}{\partial x^a}}p_b p_c.$$

The question of the integrability of the geodesic equation has been
investigated by several in the authors context of Stationary
Axially-symmetric Vacuum (SAV) space-time in \cite{Chaz},
\cite{Curz}, \cite{DubArAlf} and \cite{SamDol}. The properties of
stationarity and axial-symmetry suggest that every SAV space- time
possesses a pair Killing vector fields $\partial_t$ and $\partial
_\phi$, and  a pair first integrals $p_t=-E$, and $p_\phi =L$,
corresponding to energy $E$  and azimuthal angular momentum $L$. The
problems of geodesic motion reduces to two dimensional Hamiltonian
system in the meridian plane $(\rho, \, z)$. The existence of
additional independent integral of motion in involution with others
integrals $(E, \, L, \, H )$ lead to Liouvillian integrability. The
interesting case Kerr space-time - a SAV space-time, describing
rotating black hole - admits two-rank irreducible Killing tensor,
that is why  four constant of motion which is a quadratic in the
momenta - the so-called Carter's constants, \cite{Carter} and
\cite{ColHugh}. The important question is when the Carter's constant
exists in other SAV- models.

Geodesic motion on Zipoy-Voorhees space-times known as a
$\gamma$-metrics, has been studied via numerical  methods, in
\cite{LucGer1} and \cite{LucGer2} and analytically in
\cite{MaPrySta} and  \cite{Vol}. The Zipoy-Voorhees metrics with a
parameter $\delta$ have as  special cases both Minkowski flat
space-time ($\delta=0$) and Schwarzschild ($\delta=1$) space-time.
Important for study of Zipoy-Voorhees metrics is that finding
obvious regularity in geodesics across large parts of the parameter
$\delta$  and the finding of approximate constants  of motion,  is
not enough  for  proof of the Liouville integrability, i.e. for the
existence of fourth constant of motion.

In  \cite{Sota} it was shown numerically that the Chazy-Curzon and
Zipoy-Voorhees space-time exhibit obviously-regular dynamics in
special  parameter $\delta $ values. After further research in
\cite{Brink2}, Brink suggests that the Zipoy-Voorhees metrics may
give rise to integrable systems for $\delta\ne 0;\, 1$. In
\cite{Brink3} and  \cite{Brink4} Brink show that all SAV metrics
admits irreducible second-order Killing tensors are necessarily of
Petrov type D. This key result implies that for generic
Zipoy-Voorhees metrics without Petrov D property, do not permit in
general for any new constant of geodesic motion which is quadratic
in momentum variables.

 In 2012 in \cite{LucGer1} using numerical calculations it was shown that the Hamiltonian phase space of the Zipoy-Voorhees system has the features expected of a perturbed non-integrable system such as chaotic layers, Birkhoff chains, and ``stickiness'' in the rotation number \cite{LucGer1} and \cite{LucGer2}. Chaotic layers are mos visible near to tha Lyapunov orbit - near to the unstable periodic orbit at the border between bound and plunging orbits (see \cite{SamDol} for details).

\section{Formulation of the problem}

The Chazy-Curzon space-time model are  part of the class of Static
Axially-symmetric Vacuum (StAV) metrics - a subclass of SAV metrics.
In cylindrical coordinates $(t,\, \rho,\, z,\, \phi ) $ the StAV
line element  can be written in  Weyl form as
$$s^2=g_{\mu\nu}dx^{\mu}dx^{\nu}=-e^{2\psi}dt^2+e^{-2\psi}(e^{2\gamma}(d\rho ^2+dz^2)+\rho ^2d\phi ^2).$$
The geodesic motion is governed by a hamiltonian $H_4(x^{\mu},\,
p_{\nu})=1/2 g^{\mu \nu}p_{\mu}p_{\nu}=-\ds{\frac{1}{2}}$ with
$g^{\mu \nu}$ being the  inverse of metric. $H_4$ is a constant of
motion which we set to $-\ds{\frac{1}{2}}$ to test the bodies of
unit mass. The momentum variables are given by $p_{\mu}=g^{\mu
\nu}\dot{x}^{\nu}$, where the dot stands for differentiation with
respect to the proper time.

From Hamiltonian equations we have $\dot{p_{t}}=\dot{p_{\phi}}=0$
and we obtain two first integrals  $E=-p_t$ and $L=p_{\phi}$. Motion
on the meridian plane is governed by the reduced Hamiltonian $H_2$,
with the restriction, $H_2=0$ given by
\begin{equation}
\label{1.66} H_2=\Omega^2(\rho,\,
z).(\ds{\frac{1}{2}}(p^2_{\rho}+p^2_{z})+V_{eff}(\rho,\, z,\, E^2,\,
L^2)),
\end{equation}

where $\Omega = e^{(\psi-\gamma})$, $V_{eff}=-\ds{\frac{1}{2}}.\Phi$
with $\Phi = e^{-2\gamma}.(E^2-e^{2\psi}-\rho ^{-2}e^{4\psi}L^2)$.
The contour $V_{eff}=0$ defines a curve of zero velocity. The
geodesics are found by solving Hamiltonian equations
\begin{equation}
\label{1.1} \dot{\rho}=\ds{\frac{\partial H_2}{\partial
p_{\rho}}},\, \dot{z}=\ds{\frac{\partial H_2}{\partial p_{z}}},\,
\dot{p_{\rho}}=-\ds{\frac{\partial H_2}{\partial \rho}},\,
\dot{p_{z}}=-\ds{\frac{\partial H_2}{\partial z}}.
\end{equation}
with a choice of $E$ and $L$ and initial conditions for $
p_{\rho},\, p_z, \, \rho,\, z$ such that $H_2=0$.

Einstein's vacuum field equations for the line part reduce to
Laplace equation\\ $\partial_{\rho
\rho}\psi+\rho^{-1}\partial_{\rho}\psi +\partial_{zz}\psi =0$ with
$\gamma_{,\rho}=\rho (\psi ^2_{,\rho}+\psi ^2_{,z})$ and $\gamma_{,
z}=2\rho \psi_{,\rho}\psi_{, z}$.

The solution of the Laplace equation is
 $$\psi=-\ds{\frac{m}{\sqrt{\rho^2+z^2}}},\, \gamma=-\ds{\frac{m^2\rho^2}{2{(\rho^2+z^2})^4}}.$$
 This is a single-particle Chazy-Curzon solution.
 In this paper, we establish, that the Hamiltonian system, describing this motion is not Liouville- integrable .

Next we will find the points of equilibrium of the system,
describing geodesic motion of Chazy-Curzon space-time.  We take
$m=1$, by scaling, and we have
$$H_2( p_{\rho}, \, p_z,\, \rho, \, z)=\Omega^2(\rho,\, z).(\ds{\frac{1}{2}}(p^2_{\rho}+p^2_{z})+V_{eff}(\rho,\, z,\, E^2,\, L^2)),$$
 where $\Omega(\rho,\, z)=\ds{e^{ (\frac{-1}{\sqrt{\rho ^2+z^2}} +\frac{\rho ^2}{2(\rho ^2+z^2)^2})}} $ and change the time
 \begin{equation}
\label{1.11} d\tau=\Omega^2 (\rho,\, z)\,dt
\end{equation}
for the Hamiltonian  and the equations  (\ref{1.1}).

We found  the new ``Hamiltonian'' -  change of  the time variable is
not a canonical transformation, but preserves the integrability
properties (see \cite{HGDR} for details)
\begin{equation}
\label{1.12} F( p_{\rho}, \, p_z,\, \rho, \, z
)=\ds{\frac{1}{2}}(p^2_{\rho}+p^2_{z})+\frac{1}{2}-\frac{E^2}{2}e^\frac{2}{\sqrt{\rho^2+z^2}}+\frac{L^2}{2\rho^2}e^\frac{-2}{\sqrt{\rho^2+z^2}}
\end{equation}
and we obtain
$$\ds{\frac{\partial F( p_{\rho}, \, p_z,\, \rho, \, z)}{\partial p_{\rho}}}=p_{\rho},$$
$$\ds{\frac{\partial F( p_{\rho}, \, p_z,\, \rho, \, z)}{\partial p_{z}}}=p_z,$$

$$ -\ds{\frac{\partial F( p_{\rho}, \, p_z,\, \rho, \, z)}{\partial \rho}}=-\frac{E^2\rho ^2 e^\frac{2}{\sqrt{\rho^2+z^2}}}{(\rho ^2+z^2)^{3/2}}-\frac{L^2 e^\frac{-2}{\sqrt{\rho^2+z^2}}}{(\rho ^2+z^2)^{3/2}\rho}+\frac{L^2 e^\frac{-2}{\sqrt{\rho^2+z^2}}}{\rho^3} $$
and
$$-\ds{\frac{\partial F( p_{\rho}, \, p_z,\, \rho, \, z)}{\partial z}}=-z.(\frac{E^2\rho ^2 e^\frac{2}{\sqrt{\rho^2+z^2}}+L^2 e^\frac{-2}{\sqrt{\rho^2+z^2}}}{(\rho ^2+z^2)^{3/2}\rho^2}).$$
The substitution (\ref{1.11}), transforms the Hamiltonian system to
the  system ($`=\frac{d}{d\tau}$)
\begin{eqnarray}
\label{1.13}
\rho ' & = &  p_{\rho}, \nonumber \\
z ' & = &  p_{z}, \nonumber \\
p_{\rho}'  & = &  -\frac{E^2\rho ^2 e^\frac{2}{\sqrt{\rho^2+z^2}}}{(\rho ^2+z^2)^{3/2}}-\frac{L^2 e^\frac{-2}{\sqrt{\rho^2+z^2}}}{(\rho ^2+z^2)^{3/2}\rho}+\frac{L^2 e^\frac{-2}{\sqrt{\rho^2+z^2}}}{\rho^3}, \\
p_{z}' & = & -z.(\frac{E^2\rho ^2 e^\frac{2}{\sqrt{\rho^2+z^2}}+L^2
e^\frac{-2}{\sqrt{\rho^2+z^2}}}{(\rho
^2+z^2)^{3/2}\rho^2}).\nonumber
\end{eqnarray}
Then the conditions for the equilibrium are $p_{\rho}=p_z=z=0$ and
$$E^2\rho ^2e^{\frac{4}{\rho}}-L^2(\rho -1)=0,$$
 from $F=0$ we get
$$F(0,\, 0,\, \rho,\, 0)=-\frac{E^2}{2}e^{\frac{2}{\rho}}+\frac{L^2}{2\rho^2}e^{\frac{-2}{\rho}}+\frac{1}{2}=0.$$
Solving this linear system for $L^2$  and $E^2$ we obtain for the
equilibrium

$$L^2=-\frac{\rho ^2.e^{\frac{2}{\rho}}}{\rho -2},\, E^2=-\frac{(\rho -1).e^{\frac{-2}{\rho}}}{\rho -2}.$$
Finally, if the point $(p_{\rho},\, p_z,\, z, \, \rho)$ is an
equilibrium for Chazy- Curzon space-time Hamiltonian system, then
$p_{\rho}=p_z=z=0 $ and $\rho $ satisfies $L^2=-\frac{\rho
^2.e^{\frac{2}{\rho}}}{\rho -2}$  and $E^2=-\frac{(\rho
-1).e^{\frac{-2}{\rho}}}{\rho -2} $. For real $E^2$ and $L^2$ we
have $1<\rho<2$.

\section{Main result}

\begin{prop}
\label{pr1}
\label{pr1} On the manifold $P:=\{p_{\rho}=p_z=z=0,\, \rho\}$,
invariant under $X_F$ the system  (\ref{1.13}) has only  periodical
solutions around equilibrium $p_{\rho}=p_z=z=0$, $\rho_0\in(1,\, 2)$
with restrictions $L^2=-\frac{\rho_0^2
.e^{\frac{2}{\rho_0}}}{\rho_0-2}$, $E^2=-\frac{(\rho_0
-1).e^{\frac{-2}{\rho_0}}}{\rho_0 -2} $ .
\end{prop}
{\bf Proof:} For the potential  $v(\rho,\, z,\, E^2,\, L^2)$ we have
$$v(\rho,\, z,\, E^2,\, L^2)=\frac{1}{2}-\frac{E^2}{2}e^\frac{2}{\sqrt{\rho^2+z^2}}+\frac{L^2}{2\rho^2}e^\frac{-2}{\sqrt{\rho^2+z^2}}$$
around the points of equilibrium  we have $v(\rho, \, 0,
\,-\frac{(\rho -1).e^{\frac{-2}{\rho}}}{\rho -2}, \, -\frac{\rho
^2.e^{\frac{2}{\rho}}}{\rho -2})=0$. We compute  the second
derivative near the equilibrium,
$$\ds{\frac{\partial ^2 v(\rho, \, 0, \, -\frac{(\rho -1).e^{\frac{-2}{\rho}}}{\rho -2}, \, -\frac{\rho ^2.e^{\frac{2}{\rho}}}{\rho -2})}{\partial \rho ^2}}=-\ds{\frac{\rho^2-6\rho+4}{\rho ^4(\rho-2)}}.$$

We have only periodic solutions  around the equilibrium if   \\
$ v(\rho)=v(\rho, \, 0, \, -\frac{(\rho
-1).e^{\frac{-2}{\rho}}}{\rho -2}, \, -\frac{\rho
^2.e^{\frac{2}{\rho}}}{\rho -2})=0$,
$v'(\rho)=\ds{\frac{\partial  v(\rho, \, 0, \, -\frac{(\rho -1).e^{\frac{-2}{\rho}}}{\rho -2}, \, -\frac{\rho ^2.e^{\frac{2}{\rho}}}{\rho -2})}{\partial \rho }}=0$ and \\
$-v''(\rho)=-\ds{\frac{\partial ^2 v(\rho, \, 0, \, -\frac{(\rho
-1).e^{\frac{-2}{\rho}}}{\rho -2}, \, -\frac{\rho
^2.e^{\frac{2}{\rho}}}{\rho -2})}{\partial \rho
^2}}=\ds{\frac{\rho^2-6\rho+4}{\rho ^4(\rho-2)}}>0$, which is true
for  $\rho \in (1,\, 2)$.

 Important for the integrability of the system is the behavior of the period function,
\begin{equation}
\label{1.7} T=2\int_{\rho_-}^{\rho_+}\frac{d\rho}{\sqrt{0- v(\rho ,\
0, L^2,\ \, E^2)}},\nonumber
\end{equation}
where $\rho_{-}$ and $\rho_{+}$ are the roots of $0-v(\rho ,\ 0,
L^2,\ \, E^2)=0$,  $ v(\rho ,\ 0, L^2,\ \, E^2)$ is the potential
function, and $L^2=-\frac{\rho ^2.e^{\frac{2}{\rho}}}{\rho -2}$,
$E^2=-\frac{(\rho -1).e^{\frac{-2}{\rho}}}{\rho -2}$  are conditions
for equilibrium.

\begin{prop}
$T=2g(\rho_0)\log{\eta}+\Phi(\epsilon,\, \delta)$, where
$\Phi(\epsilon,\, \delta)$ is analytical function,
$\eta=\frac{\epsilon}{\delta}$, and
$g(\rho_0)=\frac{\sqrt{2}}{\sqrt{-v''(\rho_0)}}$.
\end{prop}
{\bf Proof:} Around the point of equilibrium $(  p_{\rho}, \, p_z\,
\rho,\, z)=( 0,\, 0,\, \rho_0,\, 0)$, the period is
\begin{eqnarray*}
T & = & 2\int_{\rho_-}^{\rho_+}\frac{d\rho}{\sqrt{0- v(\rho ,\ 0, \frac{\rho ^2.e^{\frac{2}{\rho}}}{\rho -2},\ \, \frac{(\rho -1).e^{\frac{-2}{\rho}}}{\rho -2})}}\\
   & = & 2\int_{\rho_0-\delta}^{\rho_0+\epsilon}\frac{d\rho}{\sqrt{0- v(\rho ,\ 0, \frac{\rho ^2.e^{\frac{2}{\rho}}}{\rho -2},\ \, \frac{(\rho -1).e^{\frac{-2}{\rho}}}{\rho -2})}}\\
  & =  & 2 g(\rho_0)\int_{\rho_0-\delta}^{\rho_0+\epsilon}\left(\frac{1}{\rho -\rho_0} \right)d\rho +\Phi(\epsilon ,\, \delta)\\
   & = & 2g(\rho_0) (\log(\epsilon)-\log(\delta))+\Phi(\epsilon ,\, \delta)\\
   & = & 2g(\rho_0)\log(\epsilon /\delta)+\Phi(\epsilon ,\, \delta)\\
   & = & 2g(\rho_0)\log{\eta}+\Phi(\epsilon,\, \delta),
\end{eqnarray*}
where $\Phi(\epsilon,\, \delta)$ is analytical function,
$\eta=\frac{\epsilon}{\delta}$,
$g(\rho_0)=\frac{\sqrt{2}}{\sqrt{-v''(\rho_0)}}$ \\ and
$-v''(\rho_0)=-\ds{\frac{\partial ^2 v(\rho_0, \, 0, \,
-\frac{(\rho_0 -1).e^{\frac{-2}{\rho_0}}}{\rho_0 -2}, \,
-\frac{\rho_0 ^2.e^{\frac{2}{\rho_0}}}{\rho_0 -2})}{\partial \rho
^2}}>0$ for $\rho_0 \in (1,\, 2)$.

The last condition  follows from Taylor's formula applied to $v(\rho
)$ around  $\rho_0$. We have

\begin{eqnarray}
v(\rho  ) & = & v(\rho_0)+v'(\rho_0)(\rho-\rho_0)\nonumber \\
            & + & \frac{1}{2}.v''(\rho_0)(\rho-\rho_0)^2+O((\rho-\rho_0)^3) \nonumber \\
            & = & \frac{1}{2}.v''(\rho_0)(\rho-\rho_0)^2+O((\rho-\rho_0)^3).\nonumber
\end{eqnarray}
because it  $v(\rho_0)=v'(\rho_0)=0$ for the points of equilibrium.
Then,
\begin{equation}
T=2g(\rho_0)\log{\eta}+\Phi(\epsilon,\, \delta),\nonumber
\end{equation}
around equilibrium point.

Our  goal is the following:
\begin{thm}
\label{th3}
 The  system with Hamiltonian (\ref{1.66})
is not  integrable by means of analytical first integral.

\end{thm}

We need an equivalent formulation to prove  Theorem \ref{th3}.
\begin{prop}
\label{pr0}

 The  system  (\ref{1.13})
has no additional analytical first integral.

\end{prop}

{\bf Proof:}
 For the proof  we know that there are only a periodical of solutions on $P$, of the Hamiltonian system (\ref{1.13}) on the hypersurface $F=0$. If $T$ is the period function then complex continuation of the manifolds $T=const$ turns out to be infinitely branched -  this excludes the existence of a nontrivial analytic integral on any open subset of the complex domain where this infinite branching is true.
We need to show that if $G$ is a smooth function on open set $U$
such that $V=U\cap (F=0)$ is $X_F$ invariant, $\{F,\, G\}=0$ and
derivative $d\{F,\, G\}=0$ on $V$, then $G$ is a function of $F$ and
$T$ on $V$.

The proof is similar like  in  the remarkable paper of J. J.
Duistermaat \cite{Duis}.

\section*{ Acknowledgements }

The author has been supported by contract No.  RD-22-719/29.03.2019
and grant  80-10-132 / 15.04.2019 from Sofia University ``St.
Kliment Ohridski''.

\end{document}